\documentclass[12pt,leqno]{article}
\usepackage{amsmath, amsthm, amsfonts, amssymb, color}
\usepackage{mathrsfs}
\theoremstyle{definition}

\newcommand{\scr}[1]{\mathscr #1}
\definecolor{wco}{rgb}{0.5,0.2,0.3}

\numberwithin{equation}{section} \theoremstyle{remark}

\newcommand{\ua}{\uparrow}

\title{{\bf Coupling and Applications}\footnote{Supported in part by WIMCS and NNSFC(10721091)}
}
\author{
{\bf Feng-Yu Wang}\\
\footnotesize{School of Mathematical Sciences,
Beijing Normal
University, Beijing 100875, China}\\
 \footnotesize{Department of Mathematics,
Swansea University, Singleton Park, SA2 8PP, UK}\\ \footnotesize{wangfy@bnu.edu.cn, F.-Y.Wang@swansea.ac.uk}
}
\begin{document}
\def\R{\mathbb R}  \def\ff{\frac} \def\ss{\sqrt} \def\BB{\mathbb
B}
\def\N{\mathbb N} \def\kk{\kappa} \def\m{{\bf m}}
\def\dd{\delta} \def\DD{\Delta} \def\vv{\varepsilon} \def\rr{\rho}
\def\<{\langle} \def\>{\rangle} \def\GG{\Gamma} \def\gg{\gamma}
  \def\nn{\nabla} \def\pp{\partial} \def\tt{\tilde}
\def\d{\text{\rm{d}}} \def\bb{\beta} \def\aa{\alpha} \def\D{\scr D}
\def\E{\mathbb E} \def\si{\sigma} \def\ess{\text{\rm{ess}}}
\def\beg{\begin} \def\beq{\begin{equation}}  \def\F{\scr F}
\def\Ric{\text{\rm{Ric}}} \def\Hess{\text{\rm{Hess}}}\def\B{\scr B}
\def\e{\text{\rm{e}}} \def\ua{\underline a} \def\OO{\Omega} \def\sE{\scr E}
\def\oo{\omega}     \def\tt{\tilde} \def\Ric{\text{\rm{Ric}}}
\def\cut{\text{\rm{cut}}} \def\P{\mathbb P} \def\ifn{I_n(f^{\bigotimes n})}
\def\C{\scr C}      \def\aaa{\mathbf{r}}     \def\r{r}
\def\gap{\text{\rm{gap}}} \def\prr{\pi_{{\bf m},\varrho}}  \def\r{\mathbf r}
\def\Z{\mathbb Z} \def\vrr{\varrho} \def\ll{\lambda}
\def\L{\scr L}\def\Tt{\tt} \def\TT{\tt}\def\II{\mathbb I}
\def\i{{\rm i}}\def\Sect{{\rm Sect}} \def\Q{\mathbb Q} \def\supp{{\rm supp}}
\def\T{\mathbf T} \def\LL{\Lambda}

\maketitle
\begin{abstract}   This paper presents a self-contained account for coupling arguments and applications in the context of  Markov processes. We first use   coupling to describe the transport problem, which leads to  the concepts of optimal coupling and probability distance (or transportation-cost),   then introduce applications of coupling to the study of ergodicity, Liouville theorem, convergence rate, gradient estimate, and Harnack inequality  for  Markov processes. \end{abstract} \noindent
 \vskip 2cm

\section{What is coupling}

A coupling for two distributions (i.e. probability measures) is nothing but a joint distribution of them. More precisely: 
\beg{defn} Let $(E,\F)$ be a measurable space, and let
$\mu,\nu\in \scr P(E)$, the set of all  probability  measures on $(E,\F)$. A probability measure $\pi$ on the product space $(E\times E,\F\times \F)$ is called a coupling of $\mu$ and $\nu$,
 if
 $$\pi(A\times E)=\mu(A),\ \ \pi(E\times A)=\nu(A),\ \ \ A\in\F.$$ \end{defn} We shall let $\C(\mu,\nu)$ to stand for the set of all couplings of $\mu$ and $\nu$. Obviously, the product measure $\mu\times\nu$ is a coupling of $\mu$ and $\nu$, which is called the independent coupling. This coupling is too simple to have broad applications, but it at least indicates the existence of coupling. Before moving to more general applications of coupling, let us present a simple example to show that even this trivial coupling could have   non-trivial applications. Throughout the paper, we shall let $\mu(f)$ denote the integral of  function $f$ w.r.t. measure $\mu$.

\beg{exa}[The FKG inequality] Let $\mu$ and $\nu$ be   probability measures on $\R$, then for any two bounded increasing functions $f$ and $g$, one has
$$\mu(fg)+\nu(fg) \ge \mu(f)\nu(g) +\nu(f)\mu(g).$$\end{exa}\beg{proof} Since by the  increasing monotone properties of $f$ and $g$ one has
$$(f(x)-f(y))(g(x)-g(y))\ge 0, \ \ x,y\in \R,$$ the desired inequality follows by taking integral w.r.t. the independent coupling $\mu\times\nu$.  \end{proof}
In the remainder of this section,  we shall   first link coupling to transport problem, which leads to the notions of optimal coupling and probability distances, then introduce coupling for stochastic processes.

\subsection{Coupling and transport problem}  Let $x_1, x_2,\cdots, x_n$ be $n$  places, and consider the distribution $\mu:=\{\mu_i:\ i=1,\cdots, n\}$ of some product among these places, i.e. $\mu_i$ refers to the ratio of the product at place $x_i$. We have $\mu_i\ge 0$ and $\sum_{i=1}^n \mu_i=1;$ that is, $\mu$ is a probability measure on $E:=\{1,\cdots, n\}.$   Now, due to market demand one wishes to transport the product among these places to the  target distribution $\nu:=\{\nu_i:\ 1\le i\le n\}$, which is another probability measure on $E$. Let $\pi:=\{\pi_{ij}: 1\le i,j\le n\}$ be a transport scheme, where $\pi_{ij}$  refers to the amount to be transported from place $x_i$ to place $x_j$. Obviously, the scheme is exact to transport the product from distribution $\mu$ into distribution $\nu$ if and only if $\pi$ satisfies
$$\mu_i=\sum_{j=1}^n \pi_{ij},\ \ \nu_j =\sum_{i=1}^n \pi_{ij},\ \ \ 1\le i,j\le n.$$ Thus, a scheme  transporting from $\mu$ to $\nu$ is nothing but a coupling of $\mu$ and $\nu$, and vice versa.

Now, suppose $\rr_{ij}$ is the cost to transport a unit product from place $x_i$ to place $x_j$. Then it is reasonable that $\rr$ gives rise to a distance on $E$.
With the cost function $\rr$, the transportation cost for a  scheme $\pi$ is
$$\sum_{i,j=1}^n \rr_{ij}\pi_{ij}=\int_{E\times E} \rr\,\d\pi.$$ Therefore, the minimal transportation cost between these two distributions is
$$W_1^\rr(\mu,\nu):= \inf_{\pi\in\C(\mu,\nu)} \int_{E\times E}\rr\,\d\pi,$$ which is called the $L^1$-Wasserstein distance between $\mu$ and $\nu$ induced by the cost function $\rr$.

In  general, let $(E,\F)$ be a measurable space and let $\rr$ be a non-negative measurable function on $E\times E$. For any $p\ge 1$

\beq\label{W} W_p^\rr (\mu,\nu):= \bigg\{\inf_{\pi\in \C(\mu,\nu)} \int_{E\times E} \rr^p\d\pi\bigg\}^{1/p}\end{equation} is also called the $L^p$-Wasserstein distance (or the $L^p$ transportation cost) between probability measures $\mu$ and $\nu$ induced by the cost function $\rr$. In general, $W_p^\rr$ is not really a distance
on $\scr P(E)$, but it is a distance on $\scr P_p(E):= \{\mu\in\P: \rr\in L^p(\mu\times\mu)\}$ provided $\rr$ is a distance on $E$ (see e.g. \cite{Chen}).

It is easy to see from (\ref{W}) that any coupling provides an upper bound of the transportation cost, while the following Kontorovich dual formula enables one to find lower bound estimates.

\beg{prp}[Kontorovich dual formula] Let
  $\scr F_c=\{(f,g): f,g\in \B_b(E),  f(x)\le
g(y) +\rr(x,y)^p, \ x,y\in E\},$ where $\B_b(E)$ is the set of all bounded measurable functions on $E$. Then
 $$W_p^\rr(\mu,\nu)^p= \sup_{(f,g)\in \F_c}
\{\mu(f)-\nu(g)\}.$$   \end{prp} When $(E,\rr)$ is a metric space, $\B_b(E)$ in the definition of $\F_c$    can be replaced by a sub-class of bounded measurable functions determining probability measures (e.g. bounded Lipschitzian functions), see e.g. \cite{Ra}.

\subsection{Optimal  coupling and optimal   map}

\beg{defn} Let  $\mu,\nu\in \scr P(E)$ and $\rr\ge 0$ on $E\times E$ be fixed.  If $\pi\in \C(\mu,\nu)$ reaches the infimum in (\ref{W}), then we call it an optimal coupling for the $L^p$ transportation cost. If a measurable map $\T :E\to E$ maps $\mu$ into $\nu$ (i.e. $\nu= \mu\circ\T^{-1}$), such that
$\pi(\d x,\d y): = \mu(\d x) \dd_x(\d y)$ is an optimal coupling, where $\dd_x$ is the Dirac measure at $x$, then $\T$ is called an optimal (transport) map for  the $L^p$ transportation cost.
\end{defn}

To fix (or estimate) the Wasserstein distance, it is crucial to construct the  optimal coupling or optimal map. Below we introduce some  results on existence and construction of the optimal coupling/map.

\beg{prp} Let $(E,\rr)$ be a Polish space. Then for any $\mu,\nu\in \scr P(E)$ and any $p\ge 1$, there exists an  optimal coupling. \end{prp}
The proof is fundamental. Since it is easy to see that the class $\scr C (\mu,\nu)$ is tight, for a sequence of couplings $\{\pi_n\}_{n\ge 1}$ such that
$$\lim_{n\to \infty} \pi_n(\rr^p)= W_p^\rr(\mu,\nu)^p,$$ there is a weak convergent subsequence, whose weak limit gives an optimal coupling.

As for the optimal map, let us simply mention a result of McCann for $E=\R^d$, see \cite{V} and references within for extensions and historical remarks.

\beg{thm}[\cite{McCann}]  Let  $E=\R^d, \rr(x,y)=|x-y|,$ and $p=2$. Then for any  two absolutely continuous probability measures $\mu(\d x):= f(x)\d x$ and $ \nu(\d x):=g(x)\d x$    such that   $f>0,$    there exists a unique optimal map, which is given by $T=\nn V$
 for a convex function   $V$ solving the equation
 $$f=g(\nn V) {\rm det} \nn_{ac} \nn V$$
 in the distribution sense, where   $ \nn_{ac}$
 is the gradient for the absolutely continuous part of a distribution. \end{thm}
Finally, we introduce the Wasserstein coupling which is optimal when  $\rr$ is the discrete distance on $E$; that is, this coupling is optimal for the total variation distance.

\beg{prp}[Wasserstein coupling] Let $\rr(x,y)=1_{\{x\ne y\}}$.  We have
$$ W_p^\rr (\mu,\nu)^p= \ff 1 2
\|\mu-\nu\|_{var}:=\sup_{A\in\scr F} |\mu(A)-\nu(A)|,$$  and  the Wasserstein coupling
$$\pi(\d x,\d y):= (\mu\land\nu)(\d x)\dd_x(\d y) +
\ff{(\mu-\nu)^+(\d x)(\mu-\nu)^-(\d y) }{(\mu-\nu)^-(E)}$$  is optimal, where $(\mu-\nu)^+$ and $(\mu-\nu)^-$ are the positive  and  negative parts respectively in the Hahn
 decomposition of $\mu-\nu$, and $\mu\land \nu= \mu -(\mu-\nu)^+.$ \end{prp}

\subsection{Coupling for stochastic processes} \beg{defn} Let  $X:=\{X_t\}_{t\ge 0}$ and $Y:=\{Y_t\}_{t\ge 0}$ be two stochastic processes on $E$.  A stochastic process $(\tt X, \tt Y)$ on $E\times E$ is called a coupling of them if the distributions of
$\tt X$ and $\tt Y$ coincide with those of $X$ and $Y$ respectively.  \end{defn}

Let us  observe that a coupling of two stochastic processes corresponds to a coupling of their distributions, so that the notion goes back to coupling of probability measures introduced above.

Let   $\mu$ and $\nu$ be the distributions of  $X$ and $Y$  respectively, which are probability measures on the path space
  $$W:=E^{[0,\infty)}, \ \text{equipped\ with\ product}\  \si\text{-algebra}\  \F(W):=\si\big(w\mapsto w_t:\ t\in [0,\infty)\big).$$
For any $\pi\in\scr C(\mu,\nu), (W\times W, \F(W)\times \F(W), \pi)$  is a probability space  under which
   $$(\tt X,\tt Y)(w):= (w^1, w^2),\ \ w=(w^1, w^2)\in W\times W$$   is a coupling for   $X$ and $Y$.  Conversely, the distribution of a coupling for  $X$ and $Y$  also provides a coupling for  $\mu$ and $\nu$.

\section{Some general results for Markov processes}

   Let   $P_t$  and $P_t(x,\d
  y)$ be the semigroup and transition probability kernel for a strong Markov process on a Polish space
  $E$. If   $X:=(X_t)_{t\ge 0}$  and   $Y:=(Y_t)_{t\ge 0}$  are two processes with the same transition probability kernel  $P_t(x,\d y)$, then
 $(X,Y)=(X_t,Y_t)_{t\ge 0}$ is called a coupling of the strong Markov process with coupling time
$$T_{x,y}:= \inf\{t\ge 0: X_t=Y_t\}. $$ The coupling is called  successful  if
   $T_{x,y}<\infty$ a.s. For any $\mu\in \scr P(E),$ let $\P^\mu$ be the distribution of the Markov process with initial distribution $\mu$, and let $\mu P_t$ be the marginal distribution of $\P^\mu$ at time $t$.

\beg{defn}
  If for any  $x,y\in E$, there exists a successful coupling
  starting from  $(x,y)$, then the strong Markov process is
  said to have  successful coupling  (or to have the   coupling
  property).\end{defn}

Let $$\scr T=\bigcap_{t>0} \si(
\oo\mapsto \oo_s: \ s\ge t)$$ be the tail $\si$-filed.
The following result includes some equivalent assertions for the coupling property.

\beg{thm}[\cite{CG, Lindvall, Thorisson}] Each of the
following is equivalent to the coupling property: \beg{enumerate}
\item[$(1)$] For any   $\mu,\nu\in \scr P(E),\ \lim_{t\to\infty}
\|\mu P_t -  \nu P_t\|_{var}=0.$
\item[$(2)$] All bounded time-space harmonic functions are constant, i.e. a bounded measurable function  $u$  on
 $[0,\infty)\times E$  has to be constant if
$$ u(t,\cdot)= P_s u(t+s,\cdot),\ \ s,t\ge 0.$$
 \item[$(3)$] The tail  $\si$-algebra of  is trivial,
 i.e.  $P^\mu (X\in A)=0$ or $1$ holds   for $\mu\in \scr P(E)$ and
 $A\in \scr T.$
\item[$(4)$] For any $\mu,\nu\in\scr P(E),\ \P^\mu=\P^\nu$ holds on $\scr T$.\end{enumerate}
\end{thm}

A weaker notion than the coupling property is the shift-coupling property.

\beg{defn} The strong Markov process is said to have the shift-coupling property, if for any $x,y\in E$ there is a coupling $(X,Y)$ starting at $(x,y)$ such that
$X_{T_1}=Y_{T_2}$ holds for some finite stopping times $T_1$ and $T_2$. \end{defn}

Let $$\scr I:=\big\{A\in \F(W): w\in A \ \text{implies}\ w(t+\cdot)\in A, t\ge 0\big\}$$ be the shift-invariant $\si$-field. Below are some equivalent statements for the sift-coupling property.

\beg{thm}[\cite{AT,CG,Thorisson}] Each of the
following is equivalent to the shift-coupling property: \beg{enumerate}
\item[$(5)$] For any   $\mu,\nu\in \scr P(E),\ \lim_{t\to\infty}
\ff 1 t \int_0^t\|\mu P_s -  \nu P_s\|_{var}\d s=0.$
\item[$(6)$] All bounded   harmonic functions are constant, i.e. a bounded measurable function  $f$  on
 $E$  has to be constant if
$ P_t f= f$ holds for all $t\ge 0.$
 \item[$(7)$] The invariant $\si$-algebra of the process is trivial,
 i.e.  $P^\mu (X\in A)=0$ or $1$ holds   for $\mu\in \scr P(E)$ and $A\in \scr I.$

\item[$(8)$] For any $\mu,\nu\in\scr P(E),\ \P^\mu=\P^\nu$ holds on $\scr I$.\end{enumerate}
\end{thm}

According to \cite[Theorem 5]{CW}, the coupling property and the  shift-coupling property are equivalent, and thus all above statements  (1)-(8) are equivalent, provided there exist $s,t>0$ and increasing function $\Phi\in C([0,1])$ with $\Phi(0)<1$ such that
$$P_t f\le \Phi(P_{t+s}f),\ \ 0\le f\le 1$$ holds, where ${\rm osc}(f):= \sup f-\inf f. $

By the strong Markov property, for a coupling $(X,Y)$ with coupling time $T$, we may let $X_t=Y_t$ for $t\ge T$
 without changing the transition probability kernel; that is, letting $$\tt Y_t=\beg{cases} Y_t, &\text{if}\ t\le T,\\
 X_t, &\text{if}\ t>T,\end{cases}$$ the process $(X, \tt Y)$ is again a coupling. Therefore, for any $x,y\in E$ and any coupling  $(X,Y)$ starting at $(x,y)$ with coupling times $T_{x,y}$,  we have

\beq\label{CC}|P_t f(x)-P_t f(y)| =|\E (f(X_t)-f(\tt Y_t))| \le {\rm osc}(f) \P(T_{x,y}>t),\ \  f\in \B_b(E).\end{equation}
This implies the following assertions, which are fundamentally crucial for applications of coupling in the study of Markov processes.

\beg{enumerate} \item[(i)] If $\lim_{y\to x} \P(T_{x,y}>t)=0, x\in E,$ then $P_t$ is strong Feller, i.e.
 $P_t \B_b(E)\subset C_b(E).$ \item[(ii)] Let $\mu$ be an invariant probability measure. If the coupling time $T_{x,y}$ is measurable in $(x,y)$, then
 $$\|\nu P_t - \mu\|_{var} \le 2 \int_{E\times E} \P(T_{x,y}>t)\pi(\d x,\d y),\ \ \pi\in \C(\mu,\nu)$$ holds for  $\nu\in \scr P(E).$
  \item[(iii)] The gradient estimate
 $$|\nn P_t f(x)|:= \limsup_{y\to x} \ff{|P_t f(y)-P_tf(x)|}{\rr(x,y)} \le {\rm osc}(f)\limsup_{y\to x} \ff{\P(T_{x,y}>t)}{\rr(x,y)},\ \ x\in E$$ holds.\end{enumerate}
By constructing coupling such that $\P(T_{x,y}>t)\le C\e^{-\ll t}$ holds for some $C,\ll>0$, we derive lower bound estimate of the spectral gap in the symmetric case (see \cite{CW97a, CW97b}).

\section{Derivative formula and Harnack inequality for diffusion semigroups}

To make our argument easy to follow, we shall only consider the Brownian with drift on $\R^d$. But the main idea works well for more general SDEs,   SPDEs and Neumann semigroup on manifolds with (non-convex) boundary (see \cite{ATW, DRW, ES, L, LW, RW, W07, W10b, W10, WX10,Z} and references within).

Consider the diffusion semigroup generated by  $L:=\ff 1 2 \DD+Z\cdot \nn$  on  $\R^d$  for some
  $Z\in C^1_b(\R^d,\R^d)$.  Let   $v, x\in\R^d, t>0$  be fixed. Consider  $\nn_vP_tf(x)$, the derivative of $P_t f$ at point $x$  along direction $v$,   for   $f\in \B_b(\R^d)$.   It is well known that the diffusion process starting at $x$ can be constructed by solving the It\^o SDE
$$\d X_s= \d B_s + Z(X_s)\d s,\ \ X_0=x,$$ where $B_s$ is the $d$-dimensional Brownian motion. We have  $P_t f(x)=\E f(X_t).$

  \beg{thm}[Derivative formula]  For any   $f\in \B_b(\R^d)$ and $x,v\in \R^d$,  $$ \nn_v P_t f(x)= \ff 1 t \E \bigg\{f(X_t)\int_0^t \<(t-s)\nn_vZ(X_s) +v, \d B_s\>\bigg\},\ \ t>0.$$\end{thm}

\beg{proof}
For  any $\vv>0$,  let  $X_s^\vv$  solve the equation
  $$\d X_s^\vv= \d B_s +Z(X_s)\d s  -\ff\vv t v\,\d s,\ \ X_0^\vv=x+\vv v.$$
Then  $ X_s^\vv-X_s= \ff{\vv (t-s)}t v.$  In particular,  $X_t^\vv=X_t$.   To formulate  $P_tf(x+\vv v)$  using  $X_t^\vv$, let
  $$ \tt B_s=B_s+\int_0^s \Big\{ Z(X_r)-Z(X_r^\vv) -\ff \vv t v\Big\}\d r,\ \ s\le t,$$  which is Brownian motion  under the probability measure
  $\d\P_\vv:= R_\vv\d\P$, where
  $$R_\vv:= \exp\bigg[\int_0^t\big\<Z(X_s^\vv)-Z(X_s)+\ff \vv t v,\d B_s\big\>-\ff 1 2 \int_0^t\Big|Z(X_s^\vv)-Z(X_s)+\ff \vv t v\Big|^2\d s\bigg].$$
 Reformulate the equation of  $X_s^\vv$ using $\tt B_s$:
 $$\d X_s^\vv +\d \tt B_s +Z(X_s^\vv)\d s,\ \ X_0^\vv =x+\vv v.$$  We have
$$P_t f(x+\vv v)= \E_{\P_\vv} f(X_t^\vv)= \E[R_\vv f(X_t)].$$
 Therefore,

  \beg{equation*}\beg{split}
& \nn_v P_t f(x)= \lim_{\vv\to 0} \ff{P_tf(x+\vv v)-P_t f(x)}\vv \\
&=\E \Big\{f(X_t)\lim_{\vv\to 0}\ff{R_\vv-1}\vv\Big\}\\
 &=\ff 1 t  \E \bigg\{f(X_t)\int_0^t \<(t-s)\nn_vZ(X_s) +v, \d B_s\>\bigg\}.\end{split}\end{equation*}  \end{proof}
We remark that this kind of integration by parts  formula is known as Bismut (or Bismut-Elworthy-Li) formula. But our formulation is slightly different from the Bismut-Elworthy-Li ones using   derivative processes (see \cite{B, EL}).

 \

Next, we turn to consider the Harnack inequality of $P_t$, which enables one to compare values of  $P_t f$  at different points for  $ f>0$. To this end, one may try to ask for an inequality like
 $$P_t f(x)\le C(t,x,y) P_tf(y),\ \ x,y\in\R^d, t>0,$$  where   $C: (0,\infty)\times\R^{2d}\to (0,\infty)$  is independent of  $f$.
It turns out that this inequality  is too strong to be true even for $Z=0$ (see \cite{W06} for an criterion on  existence of this inequality) . Therefore, people wish to establish weaker versions of the Harnack inequality. Using maximum principle  Li-Yau \cite{LY} established their dimension-dependent Harnack inequality with a time-shift, while using a gradient estimate argument the author \cite{W97} found a dimension-free Harnack inequality with   powers. Both inequalities have been widely applied in the study of heat kernel estimates, functional/cost inequalities and  contractivity properties of diffusion semigroups, but the latter applies also to infinite dimensional models, see
\cite{AZ, ATW, DRW, K, W06b, W07, W10b, W10, WX10} and references within.
Below, we shall introduce a coupling method for the dimension-free Harnack inequality.

   Let  $\eta$  be a positive continuous function.  Consider the coupling
  \beg{equation*}\beg{split} &\d X_s = Z(X_s)\d s +  \d B_s,\ \ X_0=x,\\
&\d Y_s= \bigg(b(Y_s) + \eta_s \cdot \ff{X_s-Y_s}{|X_s-Y_s|}\bigg)\d s +\d B_s,\ \ Y_0=y.\end{split}\end{equation*}   The additional drift $\eta_s \cdot \ff{X_s-Y_s}{|X_s-Y_s|}$ in the second equation forces $Y_t$ moves to $X_t$, and with a proper choice of function $\eta$, the force will be strong enough to make the two process move together before time $t$. We shall solve the second equation up to the coupling time
  $$\tau:=\inf\{s\ge 0:\ X_s=Y_s\}$$   and let   $X_s=Y_s$  for  $s\ge \tau$.  Assume that
\beq\label{K} \<Z(x_1)-Z(x_2),x_1-x_2\>\le K |x_1-x_2|^2,\ \ x_1,x_2\in\R^d\end{equation}  holds for some constant  $K$. Then
 $$\d |X_s-Y_s|\le \big\{K|X_s-Y_s| -\eta_s\big\}\d s,\ \ s\le\tau.$$ This implies that
 $$\e^{-K(\tau\land t)} |X_{t\land \tau} -Y_{t\land \tau}|\le |x-y| -\int_0^{t\land \tau } \e^{-Ks}\eta_s\d s.$$  Taking
 $$ \eta_s= \ff{|x-y| \e^{-Ks}}{\int_0^t\e^{-2Ks}\d s},\ s\ge 0,$$
we see   $|x-y|-\int_0^t\e^{-Ks}\eta_s\d s=0$,  so that   $\tau\le t$.  Now, let
 $$R= \exp\bigg[-\int_0^{\tau}\ff{ \eta_s}{|X_s-Y_s|}\<X_s-Y_s, \d B_s\>-\ff 1 2 \int_0^\tau \eta_s^2 \d s\bigg].$$
By the Girsanov theorem,  under the probability  $R\d\P$, the process  $Y_t$  is associated to  $P_t$.  Therefore,
$$P_t f (y)= \E [R f(Y_t)]=\E [Rf(X_t)]
 \le \big(P_tf^p(x)\big)^{1/p}\big(\E R^{p/(p-1)}\big)^{(p-1)/p}.$$
By estimating  moments of  $R$, we prove  the following result.

\beg{thm}[Dimension-free Harnack inequality]  If $(\ref{K})$ holds for some constant $K\in\R$, then \beq\label{H}\big(P_t f(x)\big)^p\le \big(P_t f^p(y)\big)\exp\bigg[\ff{pK|x-y|^2}{2(p-1)(1-\e^{-Kt})}\bigg]\end{equation}
 holds  for  $p>1$, non-negative function $f$  and   $x,y\in \R^d, t>0$.\end{thm}

 According to \cite{W10b},  for any $p>1$ the Harnack inequality (\ref{H}) implies the log-Harack inequality
$$P_t\log f(x)\le \log P_t f(y) +\ff{K|x-y|^2}{2(1-\e^{-Kt})},\ \ \ x,y\in\R^d,f\in \B_b(\R^d), f\ge 1.$$ Below, we present a simple extension of this inequality to the case with a non-constant diffusion coefficient.

\beg{thm}[\cite{RW,W10}] \label{TT} Let $\si: \R^d\to \R^d\otimes\R^d$ be   Lipschitzian such that $\si^*\si\ge \ll I$ and
\beq\label{AA}\|\si(x)-\si(y)\|_{HS}^2 + 2\<x-y, Z(x)-Z(y)\>\le K|x-y|^2,\ \ x,y\in \R^d\end{equation} hold for some constants $\ll>0$ and $K\in\R$. Then the semigroup $P_t$ generated by
$$L:=\ff 1 2 \sum_{i,j=1}^d (\si^*\si)_{ij}\pp_i\pp_j +\sum_{i=1}^d Z_i \pp_i$$ satisfies the log-Harnack inequality
$$P_t\log f(x)\le \log P_t f(y) +\ff{K|x-y|^2}{2\ll(1-\e^{-Kt})},\ \ \ x,y\in\R^d,f\in \B_b(\R^d), f\ge 1.$$\end{thm}

There are two different ways to prove this result using coupling, one is due to \cite{RW} through an $L^2$-gradient  estimate, the other is due to \cite{W10}
using coupling and Girsanov theorem. Let us briefly introduce the main ideas of these two arguments  respectively.

\beg{proof}[Proof of Theorem \ref{TT} using gradient estimate] Consider the coupling
\beg{equation*}\beg{split} &\d X_t= Z(X_t)\d t+\si(X_t)\d B_t,\ \ X_0=x,\\
&\d Y_t=Z(Y_t)\d t+\si(Y_t)\d B_t,\ \ Y_0=y.\end{split}\end{equation*} It follows from the It\^o formula and (\ref{AA}) that
$$\E |X_t-Y_t|^2\le \e^{K|x-y|^2} |x-y|.$$ Combining this with the Schwartz inequality we obtain  the $L^2$-gradient estimate
$$|\nn P_t f(x)|^2 =\lim_{y\to x} \Big(\ff{|\E(f(X_t)-f(Y_t)|}{|x-y|}\Big)^2 \le \e^{Kt} \lim_{y\to x}\E\ff{|f(X_t)-f(Y_t)|^2}{|X_t-Y_t|^2}
=\e^{Kt} P_t|\nn f|^2(x)$$ for $f\in C_b^1(\R^d)$. Up to an approximation argument, this implies that for $f\in C_b(\R^d)$ with $f\ge 1$, and for
$h\in C^1([0,t])$ such that $h_0=0, h_t=1$,
\beg{equation}\label{KK}\beg{split}& \ff{\d}{\d s} P_s\log P_{t-s}f(y+(x-y)h_s)\\
& = \big\{h'(s) \<\nn P_s \log P_{t-s} f, x-y\> -\ll P_s |\nn \log P_{t-s}f|^2\big\}((x-y)h_s+y)\\
&\le \big\{|h_s'|\cdot|x-y| \e^{Ks/2} P_s|\nn\log P_{t-s} f| -\ll P_s |\nn \log P_{t-s}f|^2\big\}((x-y)h_s+y)\\
&\le \ff {h_s'|^2 \e^{Ks}}{4\ll}|x-y|^2,\ \ s\in [0,t].\end{split}\end{equation} Taking
$$h_s= \ff{1-\e^{-Ks}}{1-\e^{-Kt}},\ \ s\ge 0$$ and integrating both sides of (\ref{KK}) over $[0,t]$, we prove the desired log-Harnack inequality.
\end{proof}

\beg{proof}[Proof of Theorem \ref{TT} using Girsanov theorem] Let $\xi_s= \ff 1 K (1-\e^{K(s-t)}), s\in [0,t]$. Consider the coupling
\beg{equation*}\beg{split} &\d X_s= Z(X_s)\d s+\si(X_s)\d B_s,\ \ X_0=x,\\
&\d Y_s=Z(Y_s)\d s+\si(Y_s)\d B_s+\ff{1_{[0,t)}(s)}{\xi_s} \si(Y_s)\si(X_s)^{-1}(X_s-Y_s)\d s,\ \ Y_0=y.\end{split}\end{equation*} From the assumption it is easy to see that
$$R_s:= \exp\bigg[-\int_0^s\ff 1 r \< \si(X_r)^{-1}(X_r-Y_r), \d B_r\>-\ff 1 2 \int_0^s \ff{|\si(X_r)^{-1}(X_r-Y_r)|^2}{\xi_r^2}\d r\bigg],\ \ s\in [0,t]$$ is a uniformly integrable martingale with
\beq\label{RR} \E [R_t\log R_t] \le \ff {K|x-y|^2} {2\ll (1-\e^{-Kt})}.\end{equation} Moreover, $X_t=Y_t$ holds $(R_t\,\d\P)$-a.s. Therefore, by the Girsanov theorem,
 (\ref{RR}) and the Young inequality, we obtain

\beg{equation*}\beg{split} P_t \log f(y)&= \E [R_t\log f(Y_t)]= \E [R_t \log f(X_t)]\le \log \E f(X_t) +\E [R_t\log R_t]\\
&\le \log P_t f(x)+ \ff{K|x-y|^2}{2\ll(1-\e^{-Kt})}.\end{split}\end{equation*}
\end{proof}

\section{Coupling  for jump processes and applications}

For a jump process, the path will be essentially
changed if a non-trivial absolutely continuous drift   is added. This means that the coupling we constructed above for diffusions with an additional drift  is no longer valid in the jump case.    Intuitively, what we can do is to add a
 $``$random jump"  in stead of  a drift. This leads to the study of   

\subsection{Quasi-invariance of random shifts} Let $X$  be a
jump process on  $\R^d$, let   $\xi$ be a random
variable on   $\R^d$,  and  let  $\tau$ be  a random time. We aim to find conditions to ensure that
the distribution of $X+\xi 1_{[\tau,\infty)}$ is absolutely
continuous with respect to that of   $X$.

We start from   a very simple jump process, i.e.  the
compound L\'evy process. $L^0$ be the  compound Poisson process on
 $\R^d$  with L\'evy measure  $\nu_0$. Let $\LL_0$ be the distribution of
  $L_0$, which is a probability measure on the path space
$$W:= \Big\{\sum_{i=1}^\infty x_i1_{[t_i,\infty)}: x_i\in
\R^d\setminus\{0\}, \ 0\le t_i\uparrow \infty\ \text{as}
 \ i\uparrow \infty\Big\}.$$  Let
  $\DD\oo_t = \oo_{t}-\oo_{t-}$ for   $\oo\in W$  and   $t>0$.

\beg{thm}[\cite{W11}] \label{Tjump} The distribution of   $L_0+\xi 1_{[\tau,\infty)}$ is
absolutely continuous with respect to  $\LL_0$ if and only if
the joint distribution of  $(L_0, \xi,\tau)$ has the form
$$\LL_0(\d\oo) \dd_0(\d z) \Theta(\oo,\d t)+ g(\oo, z,
t)\LL_0(\d\oo)\nu(\d z)\d t,$$ where  $g$  is a non-negative
measurable function on   $W\times \R^d\times [0,\infty)$, and
  $\Theta(\oo,\d t)$  is a transition measure from   $W$
to   $[0,\infty)$. In this case, the distribution of
$L_0+\xi 1_{[\tau,\infty)}$  is
\beq\label{MM} \bigg\{\P(\xi=0) + \sum_{\DD \oo_t\ne 0}
g\big(\oo-\DD\oo_t1_{[t,\infty)}, \DD\oo_t,
t\big)\bigg\}\LL_0(\d\oo).\end{equation}   \end{thm}

We note that (\ref{MM}) is an revision of  the Mecke formula on Poisson spaces.
By using quasi-invariant random shifts given in Theorem \ref{Tjump}, we are able to investigate

\subsection{ Coupling property for  O-U processes with jump}
Let  $L:=\{L_t\}_{t\ge 0}$ be the  L\'evy
process with L\'evy measure   $\nu$
(Possibly also with Gaussian and drift parts).  Let
 $ A$ be  a  $d\times d$-matrix.
 Let  $P_t$  and  $P_t(x,\d y)$  be the transition semigroup and transition probability kernel for
 the solution to the linear SDE
 $$\d X_t= A X_t\d t +\d L_t.$$
 \beg{thm}[\cite{W11} ]
Let    $\<Ax, x\>\le 0$  hold for   $x\in \R^d$. If
 $\nu\ge \rr_0(z)\d z$ such that
$$\int_{\{|z-z_0|\le\vv\}}\rr_0(z)^{-1}\d z<\infty$$ holds for some
  $z_0\in\R^d$ and some   $\vv>0$, then  $$ \|P_t
(x,\cdot)- P_t (y,\cdot)\|_{var} \le \ff{C(1+|x-y|)}{\ss t},\ \
x,y\in \R^d, t>0$$  holds for some constant  $C>0$.
\end{thm}

\paragraph{Remark.} (a) The condition $\int_{\{|z-z_0|\le\vv\}}\rr_0(z)^{-1}\d
  z<\infty$ is  very weak, as  it holds provided  $\rr_0$  has a continuous point
  $z_0\in \R^d$  such that
  $\rr_0(z_0)>0$. Successful couplings have also been constructed in \cite{SW}  under a slightly different condition.

(b) The convergence rate we derived is sharp. To see this, let   $\nu(|\cdot|^3+1)<\infty$.
 For the compound Poisson process there exists   $c>0$  such that
 $$\|P_t(x,\cdot)-P_t(y,\cdot)\|_{var} \ge \ff{c}{\ss t},\ \ \ t\ge 1 +|x-y|^2.$$

(c) The appearance of  $1$  in the upper
 bound is essential if   $\ll:=\nu(\R^d)<\infty$, as in this
 case with probability   $\e^{-\ll t}$ the process does not
 jump before time   $t$, so   that
  $$\|P_t(x,\cdot)-P_t(y,\cdot)\|_{var}\ge 2 \e^{-\ll t},\ \
 t>0, x\ne y.$$

Similarly to  what we did for  the diffusion case, we can use the coupling argument to investigate

\subsection{Derivative formula and gradient estimate}
 Let   $\nu\ge \nu_0:=\rr_0(z)\d z$  such that   $$\ll_0:=\nu_0(\R^d)<\infty.$$  The
compound Poisson process   $L^0$  with L\'evy measure
$\nu_0$  can be formulated as
$$ L_t^0= \sum_{i=1}^{N_t} \xi_i,\ \ \ t\ge 0,$$  where
 $N_t$ is the Poisson process with rate  $\ll_0$  and
  $\{\xi_i\}$  are i.i.d. random variables independent of
  $(N_t)_{t\ge 0}$  with common distribution
$\nu_0/\ll_0.$   Let  $L^1$  be the L\'evy process independent
of   $L^0$  such that    $L:=L^0+L^1$   is the L\'evy
process with L\'evy measure  $\nu$.   Let   $\tau_i$  be the
  $i$-th jump time (or ladder time) of  $N_t$. Let   $X_t^x$  be the
  solution to the liner SDE with initial value  $x$.  Consider the
  gradient of
 $$P_t^1f(x):= \E\big\{f(X_t^x)1_{\{\tau_1\le t\}}\big\}.$$
\beg{thm}[\cite{W11b}]  Let  $\rr_0\in C^1_+(\R^d)$  such
that  $\nu(\d z)\ge \rr_0(z)\d z$  and
$$ \int_{\R^d} \sup_{|x-z|\le \vv} |\nn \rr_0|(x)\d
z<\infty$$  holds for some  $\vv>0$. Then for any
$t>0$ and   $f\in \scr B_b(\R^d),$
 $$ \nn P_t^1 f(x)= \E\Big\{f(X_t^x)1_{\{N_t\ge 1\}}\ff 1
{N_t} \sum_{i=1}^{N_t} \e^{A^*\tau_i}\nn\log \rr_0(\xi_i)\Big\}.$$
\end{thm}

Next,  by the above derivative formula and comparing the small jump part with  subordinations of the  Brownian motion,
 we obtain the following result on the gradient estimate of $P_t$, which is much stronger than the strong Feller property.

\beg{thm}[\cite{W11b}]  Let   $A\le -\theta I$ and
 $$\nu(\d z)\ge |z|^{-d} S(|z|^{-2})1_{\{|z|<r_0\}}\d z$$
hold for some   $r_0>0$ and Bernstein function   $S$
such that $S(0)=0$ and
$$\aa(t):= \int_0^\infty \ff 1 {\ss r}\e^{-tS(r)}\d
r<\infty,\ \ t>0.$$  Then there exist two constants
$c_0,c_1>0$  such that
$$\|\nn P_t f\|_\infty \le c_1 \e^{-\theta^+ t}
\aa(c_0(t\land 1))\|f\|_\infty,\ \ f\in\scr B_b(\R^d), t>0.$$
 If in
particular   $A=0$  then
 $$\|\nn P_t f\|_\infty \le c_1
\Big(\aa(c_0t)+\ff 1 {r_0}\Big)\|f\|_\infty,\ \ f\in\scr
B_b(\R^d), t>0.$$\end{thm}

Obviously, if $\lim_{r\to\infty} \ff{S(r)}{\log r}=\infty$  then $\aa(t)<\infty$ holds for all   $t>0.$   Concretely,
if   $\nu(\d z)\ge c |z|^{-(d+\aa)}1_{\{|z|<r_0\}}$  (i.e.
  $S(r)=(c r)^{\aa/2}$)  for some  $c'>0$  and
$\aa\in (0,2)$  then
 $\aa(t)\le \ff {c'}{t^{1/\aa}},\ \ t>0,$  and hence,
 $$\|\nn P_t f\|_\infty\le \ff{c'\e^{-\theta^+ t} }{(t\land 1)^{1/\aa}}\|f\|_\infty.$$
 For   detailed proofs of the above results and further developments on couplings and applications of L\'evy processes, one may check with recent papers \cite{BSW, SW0, SW, W11, W11b}.

\end{document}